\documentclass[a4paper,10pt]{article}
\usepackage[a4paper, margin=1.2in]{geometry}
\usepackage{mathrsfs,fancybox,amsfonts,latexsym,theorem,amssymb,amsmath,newlfont, epsf}
\usepackage{graphicx}
\usepackage{hyperref}
\usepackage{enumerate}
\usepackage{pgf}
\usepackage{fancyhdr}
\usepackage{times}
\usepackage{pstricks,pstricks-add,pst-math,pst-xkey} 
\usepackage{nicefrac}
\usepackage{mathtools}

\theoremstyle{plain}

\newtheorem{thm}{Theorem}[section]

\def\beginpf{\noindent {\bf Proof:} \quad}
\def\endpf{\rightline{$\square$}}

\def\CC{{\mathbb C}}
\def\NN{{\mathbb N}}

\def\RR{{\mathbb R}}

\def\<{\langle}
\def\>{\rangle}

 \DeclareMathOperator{\sech}{sech}

\title{Qualitative Description of the Particle Trajectories for the $N$-Solitons Solution of the Korteweg-de Vries Equation}
\author{Ludovick Gagnon\footnote{Sorbonnes Universit\'es, UPMC Univ Paris 06, UMR 7598, Laboratoire Jacques-Louis Lions, F-75005, Paris, France, E-mail: gagnon@ljll.math.upmc.fr. LG was supported by FQRNT and by ERC advanced grant 266907 (CPDENL) of the 7th Research Framework Programme (FP7) \newline \vspace{0.4cm} \hspace{0.3cm} {\it{Keywords}} : Korteweg-de Vries equation, particles trajectory, N-solitons solution}}

\begin{document}

\maketitle

\begin{abstract}

The qualitative properties of the particle trajectories of the $N$-solitons solution of the KdV equation are recovered from the first order velocity field by the introduction of the stream function. Numerical simulations show an accurate depth dependance of the particles trajectories for solitary waves. Failure of the free surface kinematic boundary condition for the first order type velocity field is highlighted. \\

\vspace{-0.2cm}
\hspace{-0.55cm}
{\bf 2010 Mathematical Subject Classification:} 35C08, 76B15   

\end{abstract}

\section{Introduction}

The Korteweg-de Vries (KdV) equation, 
\begin{equation}\label{kdv}
\dfrac{1}{\sqrt{gh_0}}\partial_t \eta + \partial_x \eta + \dfrac{h_0^2}{6} \partial_x^3 \eta + \dfrac{3}{2h_0} \eta \partial_x \eta =0, \quad x\in \RR, \, \,  t\in \RR,
\end{equation}
corresponds to the first order approximation of unidirectional wave solutions to the full governing equations for homogeneous, non-viscous and irrotational fluid under the shallow water regime. The shallow water regime refers to the smallness of the amplitude parameter $\epsilon:=a/h_0$ and shallowness parameter $\delta:=h_0/\lambda$ in the regime $\epsilon=\mathcal{O}(\delta^2)$. Here, $h_0$ denotes the depth of the fluid, $a$ the typical amplitude of the waves, $\lambda$ the typical length of the waves and $g$ the gravitational constant. The Cartesian coordinates $(x,y)$ are chosen such that the body of fluid is represented by the two-dimensional domain bounded by the flat bottom $y=0$ and the fluid free surface $y=h_0+\eta(x,t)$. 

\begin{figure}[!th]
\begin{center}
	\includegraphics[height=4cm]{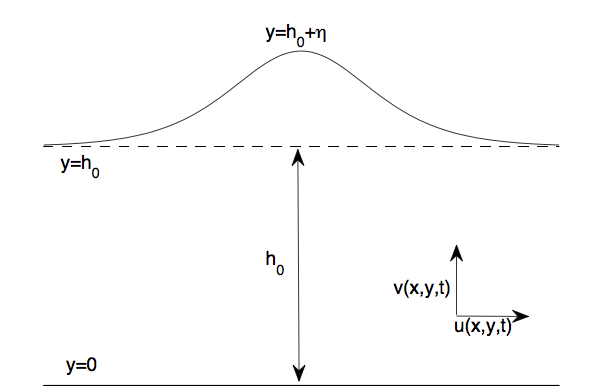}
\end{center}
\caption{\footnotesize{A soliton solution of (\ref{kdv}) represented in the Cartesian coordinates.}}
\end{figure}

The first work on the Korteweg-de Vries equation are due to John Scott Russell. Following his observation of a solitary wave on the Union Canal in 1834, a travelling wave propagating at constant speed without losing its shape, Russell published a report of his experiments and observations in 1844 \cite{soliton}. The Korteweg-de Vries equation was later derived independently by Boussinesq in 1877 \cite{Boussinesq} and by Korteweg and de Vries in 1895 \cite{KortOriginal}. 

The asymptotic behaviour of solutions of (\ref{kdv}) arising from smooth and fast decaying initial disturbance obtained by the inverse scattering transform (Gardner, Greene, Kruskal and Miura, (\cite[1967]{IST1}, \cite[1974]{IST2})), confirmed the numerical experiments of Zabusky and Kruskal (\cite[1965]{Zabusky2}): a finite number of solitons are travelling to the right while an oscillating and decaying wavetrain disperses to the left. 

A soliton of (\ref{kdv}), the approximation in the KdV regime of a solitary wave, is given by
 \begin{equation}\label{soliton}
 \eta(x,t)=a\sech^2\left(\sqrt{\dfrac{3a}{4h_0^3}}\left(x-s-\sqrt{gh_0}\left(1+\frac{a}{2h_0}\right)t \right)\right),
 \end{equation}
where $a>0$, $s\in \RR$. The following properties hold for (\ref{soliton})
\begin{enumerate}
\item The amplitude is given by $a$ and is reached at $x=s-\sqrt{gh_0}(1+\frac{a}{2h_0})t$;
\item The traveling speed is $\sqrt{gh_0}(1+\frac{a}{2h_0})t$;
\item The taller the soliton, the narrower.
\end{enumerate}

The behaviour of $N$ solitons is given by the $N$-solitons solution, representing the KdV approximation of the interactions of $N$ solitary waves propagating in the same direction. They behave like $N$ distinct solitons travelling at different constant speed except when interactions occur, in which case they are left unscathed from an interaction except for a phase shift. Moreover, if two solitons with the same order of amplitude interact, they exchange their mass and speed during the interaction.  


The Korteweg-de Vries equation not only provides an approximation of the surface elevation in shallow water but admits as well an approximation of the underlying particle trajectories. To the first order approximation in the perturbation parameter, one obtains, with the divergence-free condition, the first order velocity field of (\ref{kdv}), 
\begin{equation}
\begin{cases}\label{flow}
u(x,y,t)=\sqrt{\dfrac{g}{h_0}}\eta(x,t) \\
v(x,y,t)=-y\sqrt{\dfrac{g}{h_0}}\eta_x(x,t).
\end{cases}
\end{equation}
Since the KdV equation also holds to the second order, one recovers a second order approximation of the horizontal velocity field at the flat bottom (\cite{Whit}),
\begin{equation}\label{vitesse2}
u(x,0,t)=\sqrt{\dfrac{g}{h_0}}\left(\eta-\dfrac{1}{4h_0}\eta^2+\dfrac{h_0^2}{3}\eta_{xx}\right).
\end{equation} 
Thus, one may study the particle trajectories, for a given velocity field, associated to solutions of KdV by solving, 
\begin{equation*}
\begin{cases}
\dfrac{d X}{d t}(x,y,t)=u(X(x,y,t),Y(x,y,t),t), \quad  &t\in \RR,  \vspace{.06cm} \\
\dfrac{d Y}{d t}(x,y,t)=v(X(x,y,t),Y(x,y,t),t), \quad  &t\in \RR, \\
X(x,y,-\infty)=x, \quad \, \, \,
Y(x,y,-\infty)=y,  \quad & x\in \RR, \, \, 0\leq y \leq h_0+\eta.
\end{cases}
\end{equation*}

\begin{figure}[!th]
\begin{center}
	\includegraphics[height=8cm]{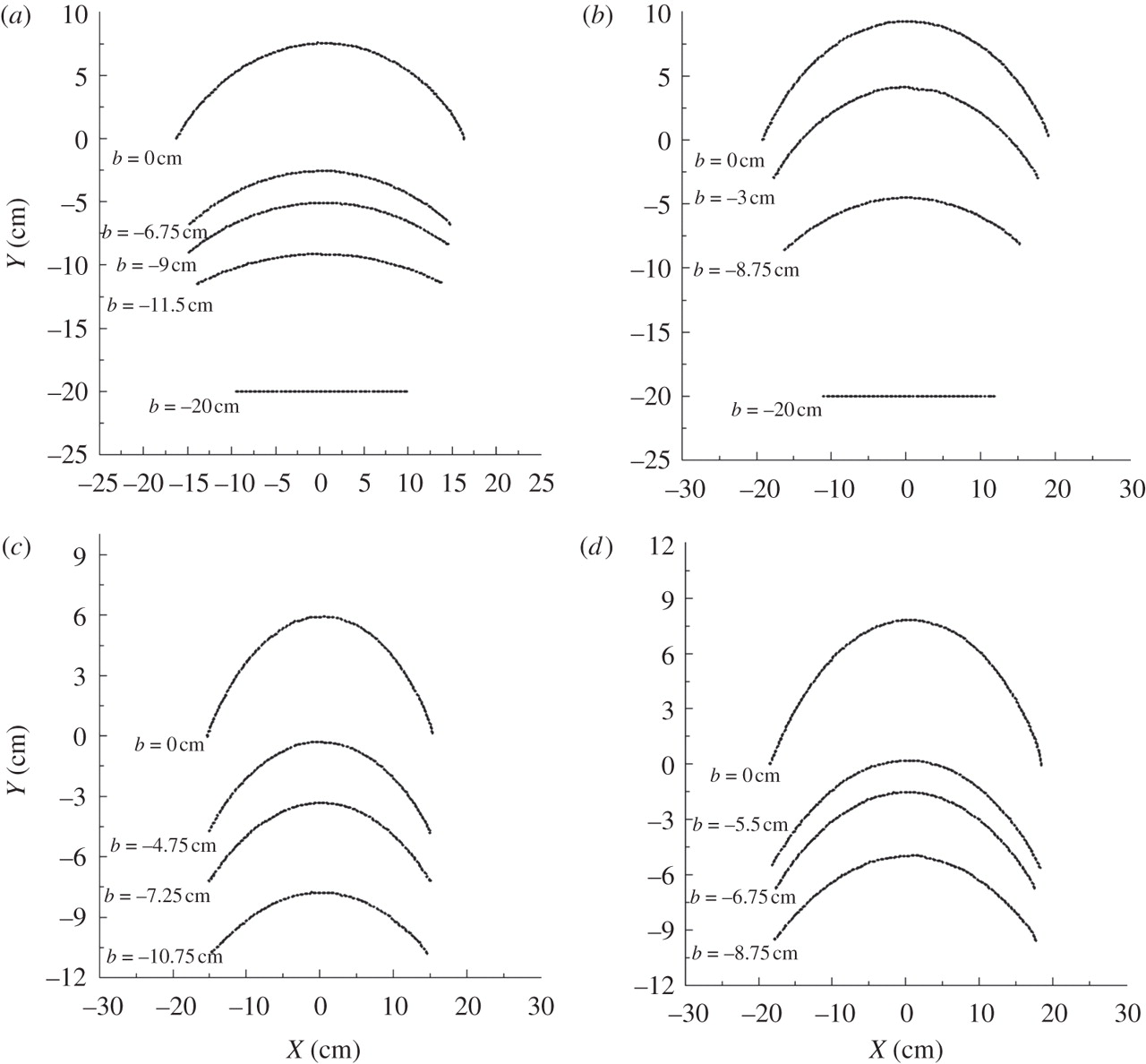}
\end{center}
\caption{\footnotesize{The orbits of water particles obtained from the experimental measurements of the polystyrene beads motions at different water levels $b$ in the four experimental wave cases (a) $h_0$=20cm, $a$=7.07cm; (b) $h_0$=20cm, $a$=8.56cm; (c) $h_0$=30cm, $a$=5.46cm; (d) $h_0$=30cm, $a$=7.56cm.}}
\label{experiment}
\end{figure}

In the case of a soliton, the particle paths obtained with the first order velocity field (\ref{flow}) fail to capture one essential property of the particle trajectories observed experimentally for solitary waves: the higher in the fluid a particle is initially located, the larger its horizontal displacement \cite{Chen}. It is a consequence of the independence on the $y$ variable of the horizontal component of the velocity field. Since the soliton solution is considered as an accurate approximation of the solitary waves solution of the full governing equations (\cite{constescher}), it is crucial to seek for higher approximations of the velocity field in order to derive better approximation of the particle trajectories. Figure \ref{experiment} illustrates the experimental results of \cite{Chen} for different solitary waves.

In the literature, two methods have been developed to recover a higher order approximation of the velocity field for the KdV equation. Introducing a stream function for the first order approximation and using the fact that the velocity potential and the stream function are harmonic conjugate, Constantin (\cite{const}) and Henry (\cite{Henry}) have obtained a higher order velocity field for a soliton and periodic solutions respectively. Qualitatives results on the particle paths were obtained in both articles. 

A careful analysis of the derivation of the KdV equation in the second order approximation allowed Borluk and Kalisch to recover from (\ref{vitesse2}) a nonlaminar velocity field underneath the surface (\cite{Kalisch}). They obtained numerically a better path of the particles in the case of the soliton solution, the cnoidals solution and the $2$-solitons solution. A similar work in the case of the Boussinesq equation had been done before in \cite{Ali,Bona}. 

The main result of this article, following the work of Constantin and Henry (\cite{const, Henry}), is the qualitative description of the particle trajectories for the $N$-solitons solution of KdV. 

\begin{thm}\label{main}
 Let $N\in \NN$ denotes the number of solitons and $a_i, \, 1\leq i \leq N$, denotes the amplitude parameter of each solitons. For sufficiently small amplitude parameters $\epsilon_i:=a_i/h_0$, the particle trajectories described by the velocity field satisfy the following:
 \begin{enumerate}[(a)]
  \item A particle on the flat bed move in straight line at a positive speed;
  \item All the particles below the surface have a positive horizontal speed;
 \end{enumerate}
 Moreover, if a particle doesn't encounter soliton interactions, we have the following: 
 \begin{enumerate}[(c)]
  \item For each such particles, there exists $t_1,\ldots,t_{2N-1}$ such that the particle move upward for,
  \[
   t\in (-\infty,t_1)\bigcup\left(\bigcup_{k=1}^{N-1}(t_{2k},t_{2k+1})\right),
  \]
 and downward for, 
 \[  
t\in \left(\bigcup_{k=1}^{N-1}(t_{2k-1},t_{2k})\right)\bigcup (t_{2N-1},\infty); 
 \]
 \end{enumerate}
  \begin{enumerate}[(d)]
 \item The particle's position at times $t_1,t_3,\ldots,t_{2N-1}$ follows the phase shift of the solitons.
 \end{enumerate}
\end{thm}

We would like to emphasize the nonlinear nature of this result. As it will be pointed out in Section 3, the phase shifts of the solitons are the result of the nonlinearity of the KdV equation. Therefore, one cannot obtain the right behaviour of the particles for the $N$-solitons solution with the superposition principle used on the 1-soliton velocity field obtained in \cite{const}. Moreover, there exists no linear approximation of solitons. This fact is of relevance when comparing the results to the periodic setting. Indeed, for the periodic travelling waves of the full governing equations, a linear approximation of the governing equations leads to nonlinear equations for the particles paths (\cite{CVillari}) and the predictions of the linear theory are reflected by what goes on in the full nonlinear system (\cite{ConstInv, ConstStrauss}). Moreover, note that as the period of the periodic travelling wave of KdV increases towards infinity, one recovers the 1-soliton (\cite{ConstBook, Johnson}).

The outline of the paper is the following. In Section 2, inspired from \cite{const,dim}, we recall the derivation of the Korteweg-de Vries to the first order in the perturbation parameter. In Section 3, we present the $N$-solitons solution and we use the method of \cite{const, Henry} in order to prove Theorem \ref{main}. This is done by noticing that this method is not limited to travelling waves. Finally, a numerical analysis of the particle trajectories is provided in Section 4. We show that the monotonicity in the $y$ variable of the horizontal displacement is recovered. Moreover, the numerical investigation indicate the possible lost of the free surface kinematic boundary condition for velocity fields of KdV.

\section{Derivation of the Korteweg-de Vries equation}

The physical assumptions to derive the Korteweg-de Vries equation are the following. Consider approximately two-dimensional waves, that is, the wave profile is the approximately the same for all cross sections taken with respect to the crest line. Assume irrotational waves in a homogeneous, non-viscous and incompressible fluid. Moreover, assume the impermeability of the bottom (bottom kinematic boundary condition), that the same particles form the surface (free surface kinematic boundary condition) and that the motion of the air and that of water is decoupled. Then, the full governing equations are given by  
\begin{equation}\label{fullgov}
\begin{cases}
u_t+uu_x+vu_y=-\frac{1}{\rho}P_x, & 0\leq y \leq h_0+\eta, \\
v_t+uv_x+vv_y=\frac{1}{\rho}P_y-g, & 0\leq y \leq h_0+\eta, \\
u_x+v_y=0, & 0\leq y \leq h_0+\eta,   \\ 
u_y-v_x=0, & 0\leq y \leq h_0+\eta,  \\ 
v=\eta_t+u\eta_x, & y=h_0+\eta, \\ 
P=P_0 & y=h_0+\eta, \\ 
v=0, & y=0.
\end{cases}
\end{equation}
Let us first set the non-dimensional full governing equations. Consider the perturbation of the pressure,
\[
 P=P_0+\rho g(h_0-y)+\rho gh_0 p,
\]
which express the change of the pressure as the waves move, and the change of variables,
\[
 x \mapsto \lambda x, \quad y \mapsto h_0 y, \quad t \mapsto \dfrac{\lambda}{\sqrt{gh_0}}t, \quad u \mapsto \sqrt{gh_0} u, \quad v  \mapsto \dfrac{\lambda}{\sqrt{gh_0^3}}v ,
\]
where $\lambda$ is the typical, or average, wavelength and where $x \mapsto \lambda x$ is to be understood as $x$ replaced by the new variable $\lambda x$. Introducing the shallowness parameter,
\[ 
\delta=h_0/\lambda,
\] 
the nondimensional full governing equations reads
\[
\begin{cases}
u_{t}+uu_{x}+vu_{y}=-p_{x}, & 0\leq {y} \leq 1+\eta/h_0, \\
\delta^2\left(v_{t}+uv_{x}+vv_{y}\right)=-p_{y}, & 0\leq {y} \leq 1+\eta/h_0, \\
u_{x}+v_{y}=0, & 0\leq {y} \leq 1+\eta/h_0,   \\ 
u_{y}-\delta^2 v_{x}=0, & 0\leq {y} \leq 1+\eta/h_0,    \\ 
v=\dfrac{1}{h_0}(\eta_{t}+u\eta_{x}), & {y}=1+\eta/h_0, \vspace{.1cm} \\ 
p=\dfrac{\eta}{h_0} & {y}=1+\eta/h_0, \\ 
v=0, & {y}=0.
\end{cases}
\]
The amplitude of the waves plays a fundamental role in the formulation of the water-wave problem. Thus, let $a$ denotes the typical, or average, amplitude of a wave and express the wave profile by $\eta \mapsto a\eta$. Then, defining the amplitude parameter by $\epsilon=a/h_0$, the nondimensional full governing equations become,
\[
\begin{cases}
u_{t}+uu_{x}+vu_{y}=-p_{x}, & 0\leq {y} \leq 1+\epsilon \eta, \\
\delta^2\left(v_{t}+uv_{x}+vv_{y}\right)=-p_{y}, & 0\leq {y} \leq 1+\epsilon \eta, \\
u_{x}+v_{y}=0, & 0\leq {y} \leq 1+\epsilon \eta,   \\ 
u_{y}-\delta^2 v_{x}=0, & 0\leq {y} \leq 1+\epsilon \eta,    \\ 
v=\epsilon(\eta_{t}+u\eta_{x}), & {y}=1+\epsilon \eta, \\ 
p=\epsilon \eta & {y}=1+\epsilon \eta, \\ 
v=0, & {y}=0.
\end{cases}
\]
Since $p$ and $v$ scale like $\epsilon$, we rescale them as well as $u$ for consistancy,
\[
 u \mapsto \epsilon u, \quad v  \mapsto \epsilon v, \quad p \mapsto \epsilon p.
\]
Then, 
\[
\begin{cases}
u_{t}+\epsilon \left( uu_{x}+vu_{y}\right)=-p_{x}, & 0\leq {y} \leq 1+\epsilon \eta, \\
\delta^2\left(v_{t}+\epsilon \left(uv_{x}+vv_{y}\right)\right)=-p_{y}, & 0\leq {y} \leq 1+\epsilon \eta, \\
u_{x}+v_{y}=0, & 0\leq {y} \leq 1+\epsilon \eta,   \\ 
u_{y}-\delta^2 v_{x}=0, & 0\leq {y} \leq 1+\epsilon \eta,    \\ 
v=\eta_{t}+\epsilon u\eta_{x}, & {y}=1+\epsilon \eta, \\ 
p=\eta & {y}=1+\epsilon \eta, \\ 
v=0, & {y}=0.
\end{cases}
\]
The regime $\epsilon=\cal{O}(\delta^2)$ arises by requiring that the above $(x,t)$-variables remain of the same order.

One removes the shallowness parameter and symmetrizes the equations by considering the change of variables
\[
 x \mapsto \dfrac{\delta}{\sqrt{\epsilon}} x, \quad  t \mapsto \dfrac{\delta}{\sqrt{\epsilon}} t, \quad v \mapsto \dfrac{1}{\delta \sqrt{\epsilon}} v,
\]
which yields,
\begin{equation}\label{fullscaled}
\begin{cases}
u_{t}+\epsilon \left( uu_{x}+vu_{y}\right)=-p_{x}, & 0\leq {y} \leq 1+\epsilon \eta, \\
v_{t}+\epsilon \left(uv_{x}+vv_{y}\right)=-p_{y}, & 0\leq {y} \leq 1+\epsilon \eta, \\
u_{x}+v_{y}=0, & 0\leq {y} \leq 1+\epsilon \eta,   \\ 
u_{y}-v_{x}=0, & 0\leq {y} \leq 1+\epsilon \eta,    \\ 
v=\eta_{t}+\epsilon u\eta_{x}, & {y}=1+\epsilon \eta, \\ 
p=\eta & {y}=1+\epsilon \eta, \\ 
v=0, & {y}=0.
\end{cases}
\end{equation}
The fourth equation of (\ref{fullscaled}) allows us to define a velocity potential $\phi(x,y,t)$ such that 
\[
u=\phi_x, \quad v=\phi_y.
\]
By restricting ourselves to waves travelling from left to right and by considering a slow time scale,
\begin{equation*}
\begin{cases}
\xi = x-t \\
T=\epsilon t
\end{cases} 
\end{equation*}
we obtain the following for $\phi$,
\begin{equation}\label{fullvar}
\begin{cases}
\epsilon \phi_{\xi \xi}+\phi_{yy}=0, & 0<y<1+\epsilon \eta, \\
\phi_y =0, & y=0, \\
\phi_y = \epsilon(-\eta_{\xi} + \epsilon \eta_T + \epsilon \phi_{\xi}\eta_{\xi}), & y=1+\epsilon \eta, \\
-\phi_{\xi} + \epsilon \phi_{T} + \dfrac{\epsilon}{2}\phi_{\xi}^2+\dfrac{\phi_y^2}{2}+\eta=0,  & y=1+\epsilon \eta. 
\end{cases} 
\end{equation}
Consider the following perturbative expansions in the amplitude parameter,
\[
 \phi(\xi,y,T) = \sum_{k=0}^\infty \epsilon^k \phi^k(\xi,y,T), \quad  \eta(\xi,T) = \sum_{k=0}^\infty \epsilon^k \eta^k(\xi,T).
\]
The zero order approximation yields, 
 \begin{equation}\label{0thorder}
\begin{cases}
\phi^0_{yy}=0, & 0<y<1, \\
\phi^0_{y} =0, & y=0, \\
\phi^0_{y} = 0,& y=1, \\
-\phi^0_{\xi}+ \dfrac{(\phi^0_{y})^2}{2}+\eta^0=0,  & y=1.
\end{cases} 
\end{equation} 
The first three equations implies, 
\[
 \phi^0_y \equiv 0, \quad 0\leq y \leq 1,
\]
while the boundary condition on $y=1$ allow us to deduce,
\begin{equation}\label{flux}
\dfrac{\partial \phi^0}{\partial \xi} (\xi,y,T)=\eta^0(\xi,T), \quad 0\leq y \leq 1.
\end{equation}
Using the first and second order of (\ref{fullvar}), one obtains that $\eta_0$ satisfies the Korteweg-de Vries equation
\[
2\eta^0_{T_0}+3\eta^0\eta^0_{\xi} +  \dfrac{1}{3}\eta^0_{\xi \xi \xi}=0.
\]
The region where the KdV balance occurs is when $\xi=O(1), T_0=O(1)$, that is, in the physical variables, when 
\[
 x\dfrac{\epsilon^{3/2}}{\lambda \delta}=O(1), \quad x-\sqrt{gh_0}t=O(1).
\]
Coming back to the physical variables, one obtains the Korteweg-de Vries equation (\ref{kdv}) and a laminar flow for solutions $\eta$ of (\ref{kdv}) 
\[
\begin{cases}
u(x,y,t)=\sqrt{\dfrac{g}{h_0}}\eta(x,t)\\ 
v(x,y,t)=-y\sqrt{\dfrac{g}{h_0}}\eta_x(x,t).
\end{cases}
\]

\section{A higher velocity field for the $N$-soliton solution}
\subsection{The $N$-solitons solution}

Let us recall the $N$-solitons solution as introduced by Hirota in \cite{Hirota} (see also \cite{Whit}) for the Korteweg-de Vries (\ref{kdv}). Let
\begin{equation}\label{expu}
\eta=-\dfrac{4h_0^3}{3}\left( \ln F\right)_{xx},
\end{equation}
where $F$ is given by the power serie expansion,
\begin{equation}\label{expension}
 F=1+\epsilon F^{(1)} + \epsilon^2 F^{(2)} + ...
\end{equation}
For $N \in \NN$, let $a_i>0, s_i \in \RR, 1\leq i \leq N$,
\[
 f_i(x,t)=\exp\left(-2\sqrt{\dfrac{3a_i}{4h_0^3}}\left(x-s_i-\sqrt{gh_0}\left(1+\frac{a_i}{2h_0}\right)t\right)\right).
\]
and
\[
 F^{(1)}=\sum_{n=1}^N f_i(x,t).
\]
Then, by replacing (\ref{expu}) in (\ref{kdv}), one obtains that only the first $N$ order term in (\ref{expension}) are nonzero. Moreover, one may express $F$ as 
\[
F=1+\sum_{n=1}^N\sum_{_{N}C_{n}}\alpha(i_1,...,i_n) \prod_{j=1}^{n} f_{i_j},
\]
with
\begin{align*}
\alpha(i_1,...,i_n)&= \prod^{(n)}_{k<l} \alpha(i_k,i_l), \qquad  \quad \textrm{ if } n\geq 2, \\
\alpha(i_k,i_l)&=\left(\dfrac{\sqrt{a_{i_k}}-\sqrt{a_{i_l}}}{\sqrt{a_{i_k}}+\sqrt{a_{i_l}}}\right)^2, \\
\alpha(i_k)&=1,  				
\end{align*}
where $\sum_{_{N}C_n}$ is the sum over all the possible combinations of $n$ indexes taken from $\{1,...,N\}$ and where 
\[
\prod^{(n)},
\] 
is the product over these $n$ indexes, provided $k<l$ when specified. It will be convenient to denote
\[
 \beta_i:=\sqrt{\dfrac{3a_i}{4h_0^3}}, \quad c_i:=\sqrt{gh_0}\left(1+\frac{a_i}{2h_0}\right).
\]

Let us explain why this solution is called the $N$-solitons solution. Consider $0<a_N<...<a_1$ and $s_N < ... < s_1 $. From (\ref{expu}), $\eta$ can be recasted as
\begin{align}\label{dec}
\eta=&4 \left( \dfrac{\sum_{n=1}^N\sum_{_{N}C_{n}} \left(\sqrt{a_{i_1}}+...+\sqrt{a_{i_n}}\right)^2 \alpha(i_1,...,i_n) \prod^{(n)} f_{i_n}}{\left(1+\sum_{n=1}^N\sum_{_{N}C_{n}}\alpha(i_1,...,i_n) \prod^{(n)} f_{i_n}\right)^2} \right. \nonumber
\\ & \quad \dfrac{+\sum_{n,m=1}^N\sum_{_{N}C_{n,m}} \left(\sqrt{a_{i_1}}+...+\sqrt{a_{i_n}}\right)\left[\left(\sqrt{a_{i_1}}+...+\sqrt{a_{i_n}}\right)   \right. }{\left(1+\sum_{n=1}^N\sum_{_{N}C_{n}}a(i_1,...,i_n) \prod^{(n)} f_{i_n}\right)^2} \nonumber \\ 
& \quad \left. \dfrac{\left. -\left(\sqrt{a_{j_1}}+...+\sqrt{a_{j_m}}\right) \right] \alpha(i_1,...,i_n)\alpha(j_1,...,j_m) \prod^{(n)}\prod^{(m)} f_{i_n} f_{j_m}}{\left(1+\sum_{n=1}^N\sum_{_{N}C_{n}}\alpha(i_1,...,i_n) \prod^{(n)} f_{i_n}\right)^2}\right).
\end{align}

First, suppose that we are in the region where $f_1\simeq 1$ and $f_k \ll 1$ for all $k\geq 2$. The behaviour of $\eta$ is given by,
 \begin{align}\label{soliton1}
\eta &\simeq 4a_1 \dfrac{f_1}{\left( 1+ f_1 \right)^2} \nonumber \\ 
&= a_1 \sech^2\left(\beta_1\left(x-s_1-c_1 t \right)\right),
\end{align}
that is, a soliton of height $a_1$ and phase $s_1$.

\begin{figure}[!th]
\begin{center}
\begin{tabular}{c}
	\includegraphics[height=6cm]{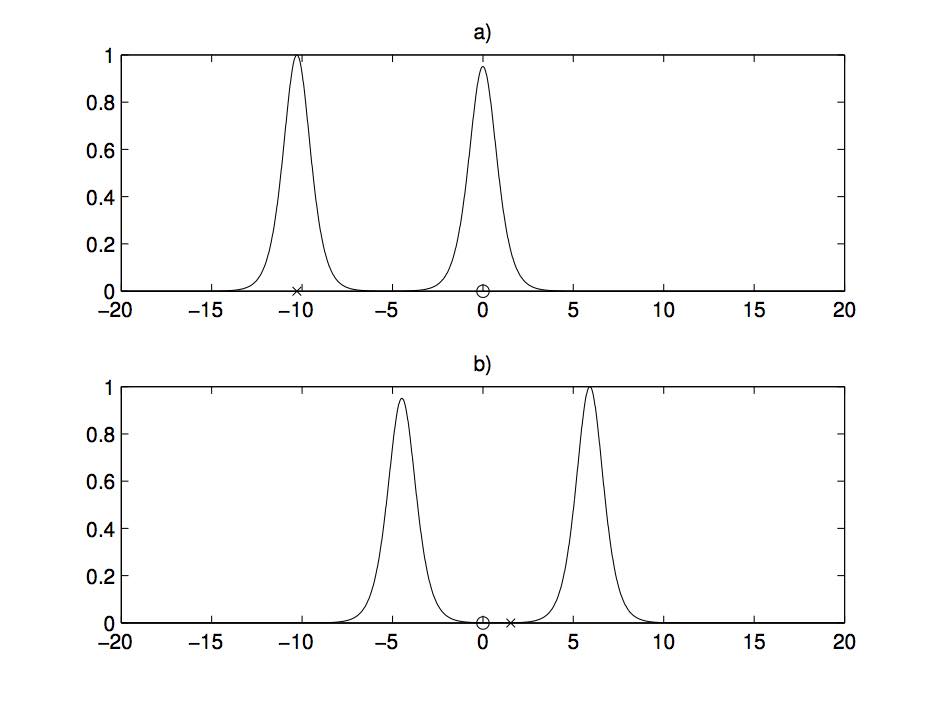}
\end{tabular}
\end{center}
\caption{\footnotesize{Interaction between two solitons. The cross (resp. circle) represents the position of the maximum of the faster (resp. slower) soliton if no interaction would have occured. Figure a) is the state of the $2$-solitons solution before the interaction and b) is the state after the interaction. The frame is fixed at the speed of the slower soliton.}}\label{interaction}
\end{figure}

Suppose now, for $2\leq k \leq N$, that we are in the region where $f_k\simeq 1 $, $f_i\gg 1$ for $1\leq i < k$ and $f_i \ll 1$ for $ k< i \leq N$. Then, 
 \begin{align}\label{shift}
\eta &\simeq 4a_k \dfrac{\alpha(1,...,k)\alpha(1,...,k-1)f_1^2\dotsm f_{k-1}^2f_k}{\left( \alpha(1,...,k-1)f_1\dotsm f_{k-1} + \alpha(1,...,k)f_1\dotsm f_{k} \right)^2} \nonumber \\ 
&= 4a_k \dfrac{\frac{\alpha(1,...,k)}{\alpha(1,...,k-1)}f_k }{\left(1+\frac{\alpha(1,...,k)}{\alpha(1,...,k-1)}f_k\right)^2}, \nonumber \\
&= a_k \sech^2\left(\beta_k\left(x-s_k-\frac{1}{2\beta_k}\ln\left(\prod^{k-1}_{i=1}\alpha(i,k)\right)-c_k t \right)\right).
\end{align}
that is, a soliton with a phase shift of $\frac{1}{2\beta_k}\ln\left(\prod^{k-1}_{i=1}\alpha(i,k)\right)$. 
From the last computations, one sees that the $N$-solitons solution behaves like $N$ distinct solitons when they are far apart from each other. Moreover, interactions between solitons must occur since the solitons are travelling at different speed. Indeed, a direct computation shows that $a_i\neq a_j$, otherwise the solution becomes a $(N-k)$-solitons solution, $k$ being the number of times $a_i=a_j, \, i<j$. Thus, up to a change of indexes in the last computations, one sees that they are left unchanged in shape after interactions, the only notable effect of the interaction being the phase shift. Figure \ref{interaction} illustrates the phase shift resulting from the interaction of 2 solitons. 

\subsection{Nonlaminar velocity field}

The velocity potential obtained in the previous section is independent of the variable $y$. Therefore, even if the soliton solution of the KdV equation provides a good approximation of the solitary wave, one cannot describe accurately the underlying structure of the particles paths.

Let us explain how one may recover a nonparallel displacement of the particles. Consider that $\phi$ is the velocity potential of \eqref{fullgov} satisfying,
\[
 \phi_x=u, \quad \phi_y=v.
\]
We introduce the stream function $\psi$ of \eqref{fullgov} satisfying
\[
\psi_y=u, \quad \psi_x=-v.
\]
For irrotational flows, the velocity potential and stream function are harmonic conjugates. Consider the complex velocity potential $f$, with $t$ acting as a parameter, such that
\[
 f(z,t):=\phi(x,y,t)+i\psi(x,y,t), \, z\in \CC,
\]
with 
\[
 f'(z,t)=u(x,y,t)-iv(x,y,t).
\]
Instead of using \eqref{0thorder} to obtain the expression an expression of the velocity potential independant of the $y$ variable, we rather solve the Laplacian equation of the velocity potential and of the stream function in the upper-half plane $\{z=x+iy\in \CC | y>0\}$ with the velocity field \eqref{flow} as a boundary condition on $y=0$. 
\begin{equation}\label{bord}
\phi_y(x,0,t)=0, \quad \phi_x(x,0,t)=\sqrt{\dfrac{g}{h_0}}\eta(x,t), \quad x\in \RR, \, \, t\in \RR
\end{equation}
If the first boundary condition is the natural non penetration condition, we use the independance of the $y$ coordinate of the horizontal component of the velocity field to obtain the second boundary condition (\cite{const, Henry}).

The boundary condition implies 
\[
 f'(x,0,t)=\sqrt{\dfrac{g}{h_0}}\eta(x,t).
\]
By the reflection principle, the analytic function $f'(z)$ has an unique extension across $y=0$. Therefore, $f$ is defined up to a constant. Using Alembert's formula, we obtain the following expressions of $u$ and $v$
\begin{align*}
 u(x,y,t)&= \dfrac{1}{2}\sqrt{\dfrac{g}{h_0}}\left(\eta(x+iy,t)+\eta(x-iy,t)\right), \\
 v(x,y,t)&= \dfrac{i}{2}\sqrt{\dfrac{g}{h_0}}\left(\eta(x+iy,t)-\eta(x-iy,t)\right).
\end{align*}
In the case where $\eta$ is the $N$-solitons solution, an explicit expression of the velocity is obtained. Indeed, let,
\begin{align*}
F(x+iy,t)=&F^c(x,y,t)-iF^s(x,y,t), \\
F(x-iy,t)=&F^c(x,y,t)+iF^s(x,y,t), 
\end{align*}
where,
\begin{align*}
F^c(x,y,t):=&1+\sum_{n=1}^N\sum_{_{N}C_{n}}\alpha(i_1,...,i_n)\left(\prod^{n}_{j=1}f_{i_j}(x,t)\right)\cos\left(\sum^{n}_{j=1}2\beta_{i_j}y\right), \\
F^s(x,y,t):=&\sum_{n=1}^N\sum_{_{N}C_{n}}\alpha(i_1,...,i_n)\left(\prod^{n}_{j=1}f_{i_j}(x,t)\right)\sin\left(\sum^{n}_{j=1}2\beta_{i_j}y\right).
\end{align*}
Then, one has,
\begin{align}
 u(x,y,t)=& \dfrac{1}{2}\dfrac{4h_0^3}{3}\sqrt{\dfrac{g}{h_0}}\left(\eta(x+iy,t)+\eta(x-iy,t)\right) \nonumber \\
 =&\dfrac{1}{2}\dfrac{4h_0^3}{3}\sqrt{\dfrac{g}{h_0}}\left(\dfrac{(F_{xx}^c-iF^s_{xx})(F^c-iF^s)-(F_{x}^c-iF^s_{x})^2}{(F^c-iF^s)^2}\right. \nonumber \\
 &\qquad \qquad \quad  +\left.  \dfrac{(F_{xx}^c+iF^s_{xx})(F^c+iF^s)-(F_{x}^c+iF^s_{x})^2}{(F^c+iF^s)^2}\right) \nonumber \\ 
 =&\textrm{Re } \left( \dfrac{4h_0^3}{3}\sqrt{\dfrac{g}{h_0}}\left(\dfrac{(F_{xx}^c-iF^s_{xx})(F^c+iF^s)[(F^c)^2+(F^s)^2]}{((F^c)^2+(F^s)^2)^2} \right. \right. \nonumber \\
 & \qquad \qquad \qquad \qquad \left.\left. -\dfrac{[(F_{x}^c-iF^s_{x})(F^c+iF^s)]^2}{((F^c)^2+(F^s)^2)^2} \right) \right) \nonumber \\ 
  =&\textrm{Re } \left( \dfrac{4h_0^3}{3}\sqrt{\dfrac{g}{h_0}}\left(\dfrac{(F_{xx}^c-iF^s_{xx})(F^c+iF^s)[(F^c)^2+(F^s)^2]}{((F^c)^2+(F^s)^2)^2} \right. \right. \nonumber \\
 & \qquad \qquad \qquad \qquad  \left.\left. -\dfrac{[(F_{x}^c-iF^s_{x})(F^c+iF^s)]^2}{((F^c)^2+(F^s)^2)^2} \right) \right) \nonumber \\ 
   =&\dfrac{4h_0^3}{3}\sqrt{\dfrac{g}{h_0}}\left(\dfrac{(F^cF_{xx}^c+F^sF^s_{xx})[(F^c)^2+(F^s)^2]}{((F^c)^2+(F^s)^2)^2} \right. \nonumber \\
 & \qquad \qquad \quad  \left. -\dfrac{[(F^cF_{x}^c+F^sF^s_{x})^2-(F^sF_x^c-F^cF^s_x)^2]}{((F^c)^2+(F^s)^2)^2} \right) \nonumber  \\ 
=&\dfrac{4h_0^3}{3}\sqrt{\dfrac{g}{h_0}}\left(\dfrac{((F^c_x)^2+F^cF_{xx}^c+(F^s_x)^2+F^sF^s_{xx})[(F^c)^2+(F^s)^2]}{((F^c)^2+(F^s)^2)^2} \right. \nonumber \\
 & \qquad \qquad \qquad \qquad\qquad  \left. -2\dfrac{(F^cF_{x}^c+F^sF^s_{x})^2}{((F^c)^2+(F^s)^2)^2} \right).  \label{u}
\end{align}

On the other hand,
\begin{align}
 v(x,y,t)=&\dfrac{i}{2}\dfrac{4h_0^3}{3}\sqrt{\dfrac{g}{h_0}}\left(\eta(x+iy,t)-\eta(x-iy,t)\right) \nonumber \\
 =&\dfrac{i}{2}\dfrac{4h_0^3}{3}\sqrt{\dfrac{g}{h_0}}\left(\dfrac{(F_{xx}^c-iF^s_{xx})(F^c-iF^s)-(F_{x}^c-iF^s_{x})^2}{(F^c-iF^s)^2}\right. \nonumber \\
&\qquad \qquad \left.  -\dfrac{(F_{xx}^c+iF^s_{xx})(F^c+iF^s)-(F_{x}^c+iF^s_{x})^2}{(F^c+iF^s)^2}\right) \nonumber \\ 
 =&-\textrm{Im } \left( \dfrac{4h_0^3}{3}\sqrt{\dfrac{g}{h_0}}\left(\dfrac{(F_{xx}^c-iF^s_{xx})(F^c+iF^s)[(F^c)^2+(F^s)^2]}{((F^c)^2+(F^s)^2)^2} \right. \right. \nonumber \\
 & \qquad \qquad \qquad \qquad \qquad  \left.\left. +\dfrac{[(F_{x}^c+iF^s_{x})(F^c-iF^s)]^2}{((F^c)^2+(F^s)^2)^2} \right) \right) \nonumber \\ 
  =& \dfrac{4h_0^3}{3}\sqrt{\dfrac{g}{h_0}}\left(\dfrac{(F^cF^s_{xx}-F^sF_{xx}^c)[(F^c)^2+(F^s)^2]}{((F^c)^2+(F^s)^2)^2} \right.  \nonumber \\
 & \qquad \qquad \qquad  \left. -2\dfrac{(F^cF_{x}^c+F^sF^s_{x})(F^cF^s_x-F^sF^c_x)}{((F^c)^2+(F^s)^2)^2} \right) \nonumber \\ 
   =&\dfrac{4h_0^3}{3}\sqrt{\dfrac{g}{h_0}}\left(\dfrac{(F^cF^c_{xy}+F^sF_{xy}^s)[(F^c)^2+(F^s)^2]}{((F^c)^2+(F^s)^2)^2}  \right. \nonumber \\
 & \qquad \qquad \quad  \left. -2\dfrac{(F^cF_{x}^c+F^sF^s_{x})(F^cF^c_y+F^sF^s_y)}{((F^c)^2+(F^s)^2)^2} \right) , \label{v}
\end{align}
the last line coming from the expression of the derivatives of $F^c$ and $F^s$, 
\begin{align*}
\dfrac{\partial^k F^c}{\partial x^k}&=\sum_{n=1}^N\sum_{_{N}C_{n}}(-2\beta_{i_1}-...-2\beta_{i_n})^k \alpha(i_1,...,i_n)\Big(\prod^{n}_{j=1}f_{i_j}(x,t)\Big)\cos\Big(\sum^{n}_{j=1}2\beta_{i_j}y\Big), \\
\dfrac{\partial^k F^s}{\partial x^k}&=\sum_{n=1}^N\sum_{_{N}C_{n}}(-2\beta_{i_1}-...-2\beta_{i_n})^k\alpha(i_1,...,i_n)\Big(\prod^{n}_{j=1}f_{i_j}(x,t)\Big)\sin\Big(\sum^{n}_{j=1}2\beta_{i_j}y\Big), \\
\dfrac{\partial F^c}{\partial y}&=\sum_{n=1}^N\sum_{_{N}C_{n}}(-2\beta_{i_1}-...-2\beta_{i_n})\alpha(i_1,...,i_n)\Big(\prod^{n}_{j=1}f_{i_j}(x,t)\Big)\sin\Big(\sum^{n}_{j=1}2\beta_{i_j}y\Big), \\
\dfrac{\partial F^s}{\partial y}&=\sum_{n=1}^N\sum_{_{N}C_{n}}(2\beta_{i_1}+...+2\beta_{i_n})\alpha(i_1,...,i_n)\Big(\prod^{n}_{j=1}f_{i_j}(x,t)\Big)\cos\Big(\sum^{n}_{j=1}2\beta_{i_j}y\Big). 
\end{align*}
\subsection{Proof of Theorem 1.1}

Before proving Theorem \ref{main}, let us introduce some notations to simplify some of the expressions. We denote 
\[
\sum:=\sum_{n=1}^N\sum_{_{N}C_{n}}, \quad \prod:=\prod^{n}_{j=1},
\]
$\sum^m$, for $m$ sums, and  $\prod^m$, for $m$ products. We shorten the expression of $\alpha(i_1,...,i_n)$ by $\alpha_{i_n}$ and define $\alpha_{i_n,j_m}:=\alpha_{i_n}\alpha_{j_m}$ and so on for product over more than two terms. Finally, we use the slightly modified Einstein notation $(\beta_{i_n})$ for a sum from $\beta_{i_1}$ to $\beta_{i_n}$, while $\beta_{i_j}$ denotes a single coefficient. Since no confusion is possible, $(\beta_{i_n})(\beta_{j_m})$ denotes the product over the two sums.

Let us turn ourselves to the proof of Theorem \ref{main}

\beginpf

(a) From (\ref{u})-(\ref{v}), we recover the original velocity field when $y=0$. 

(b) The expression of the numerator of $u$ is given by, 
\begin{align*}
 & \sum (2\beta_{i_n})^2 \alpha_{i_n} \prod f \cos((2\beta_{i_n})y) \\
 +& \sum^2 \left[3(2\beta_{i_n})^2 - (2\beta_{i_n})(2\beta_{j_m})\right]\alpha_{i_n,j_m} \prod^2 f \cos((2\beta_{i_n})y)\cos((2\beta_{j_m})y) \\
 +& \sum^3 \left[3(2\beta_{i_n})^2 - 2(2\beta_{i_n})(2\beta_{j_m})\right] \\ 
 & \times \alpha_{i_n,j_m,k_p} \prod^3 f \cos((2\beta_{i_n})y)\cos((2\beta_{j_m})y)\cos((2\beta_{k_p})y) \\
 +& \sum^4 \left[(2\beta_{i_n})^2 - (2\beta_{i_n})(2\beta_{j_m})\right] \\
 &\times \alpha_{i_n,j_m,k_p,l_q}\prod^4 f \cos((2\beta_{i_n})y)\cos((2\beta_{j_m})y)\cos((2\beta_{k_p})y)\cos((2\beta_{l_q})y) \\
 +& \sum^2 \left[(2\beta_{i_n})^2 + (2\beta_{i_n})(2\beta_{j_m})\right]\alpha_{i_n,j_m} \prod^2 f \sin((2\beta_{i_n})y)\sin((2\beta_{j_m})y)    \\  
 +& \sum^4 \left[(2\beta_{i_n})^2 - (2\beta_{i_n})(2\beta_{j_m})\right] \\ 
 & \times \alpha_{i_n,j_m,k_p,l_q} \prod^4 f \sin((2\beta_{i_n})y)\sin((2\beta_{j_m})y)\sin((2\beta_{k_p})y)\sin((2\beta_{l_q})y) \\
+& \sum^3 \left[(2\beta_{i_n})^2 + 2(2\beta_{k_p})(2\beta_{l_q})+2(2\beta_{k_p})^2-4(2\beta_{i_n})(2\beta_{k_p})\right] \\
& \times \alpha_{i_n,j_m,k_p} \prod^3 f \cos((2\beta_{i_n})y)\sin((2\beta_{j_m})y)\sin((2\beta_{k_p})y)
\end{align*}
\begin{align*}
+& \sum^4 \left[(2\beta_{i_n})^2 +(2\beta_{i_n})(2\beta_{j_m})+(2\beta_{k_p})^2 +(2\beta_{k_p})(2\beta_{l_q})-4(2\beta_{i_n})(2\beta_{k_p})\right ] \\ 
& \times  \alpha_{i_n,j_m,k_p,l_q} \prod^4 f \cos((2\beta_{i_n})y)\cos((2\beta_{j_m})y)\sin((2\beta_{k_p})y)\sin((2\beta_{l_q})y).
  \end{align*}
The coefficients of the first sum with respect to $\beta_i$ are all positive. Moreover, the coefficients appearing in the other sums are non-negative. Indeed, consider for example the coefficient in front of \\ $\cos((2\beta_{i_n})y)\cos((2\beta_{j_m})y)\sin((2\beta_{k_p})y)\sin((2\beta_{l_q})y)$ : 
\begin{align*}
& \quad \quad (2\beta_{i_n})^2 +(2\beta_{i_n})(2\beta_{j_m})+(2\beta_{k_p})^2 +(2\beta_{k_p})(2\beta_{l_q})-4(2\beta_{i_n})(2\beta_{k_p}) \\
& \quad + (2\beta_{j_m})^2 +(2\beta_{j_m})(2\beta_{i_n})+(2\beta_{l_q})^2 +(2\beta_{l_q})(2\beta_{k_p})-4(2\beta_{j_m})(2\beta_{l_q}) \\
&=((2\beta_{i_n})-(2\beta_{k_p}))^2+(2\beta_{i_n})((2\beta_{j_m})-(2\beta_{k_p}))+(2\beta_{k_p})((2\beta_{l_q})-(2\beta_{i_n})) \\
&+((2\beta_{j_m})-(2\beta_{l_q}))^2+(2\beta_{j_m})((2\beta_{i_n})-(2\beta_{i_n}))+(2\beta_{l_q})((2\beta_{k_p})-(2\beta_{j_m})) \\
&\qquad \qquad \qquad \quad =((2\beta_{i_n})+(2\beta_{j_m})-(2\beta_{k_p})-(2\beta_{l_q}))^2.
\end{align*}
Therefore, since
\[
 0\leq y \leq h_0+a, 
\]
where, 
\[
a:=\max_{\RR}\eta = \max_{1\leq i \leq N} a_i,
\]
in order to prove (b), it is sufficient to impose,  
\begin{equation}\label{hyp}
(h_0+a)\sum_{i=1}^N \beta{i} \leq \dfrac{\pi}{4}.
\end{equation}

(c) Let us give an equivalent form of $v$. Consider first
\begin{align*}
\lefteqn{ \left[F^cF^c_{xy}+F^sF^s_{xy}\right]  =4\sum  (\beta_{i_n})^2 \alpha_{i_n}\prod f \sin(2(\beta_{i_n})y)} \\ 
& \qquad \qquad  \qquad   + 4\sum^2 \left[ (\beta_{j_m})^2-(\beta_{i_n})^2\right]\alpha_{i_n,j_m}\prod^2 f \cos(2(\beta_{i_n})y)\sin(2(\beta_{j_m})y).
\end{align*}
Then coefficient associated to $(f_{i_1}\ldots f_{i_n})(f_{j_1}\ldots f_{j_m})$ in this expression is given by
\begin{align*}
4\alpha_{i_n,j_m}&\left(\left[(\beta_{j_m})^2-(\beta_{i_n})^2\right]\cos(2(\beta_{i_n})y)\sin(2(\beta_{j_m})y) \right. \\
&\left. \, + \left[(\beta_{i_n})^2-(\beta_{j_m})^2\right]\cos(2(\beta_{j_m})y)\sin(2(\beta_{i_n})y) \right),
\end{align*}
which can be written as
\[
4\alpha_{i_n,j_m}\left[(\beta_{j_m})^2-(\beta_{i_n})^2\right]\sin(2((\beta_{j_m})-(\beta_{i_n}))y).
\]
Therefore, using the symmetry of the second sum, one has 
\begin{align*}
\lefteqn{ \left[F^cF^c_{xy}+F^sF^s_{xy}\right]  =4\sum  (\beta_{i_n})^2 \alpha_{i_n}\prod f \sin(2(\beta_{i_n})y)} \\ 
& \qquad \qquad  \qquad \quad   + 2\sum^2 \left[ (\beta_{j_m})^2-(\beta_{i_n})^2\right]\alpha_{i_n,j_m}\prod^2 f \sin(2((\beta_{j_m})-(\beta_{i_n}))y), 
\end{align*}
In the same fashion, one obtains the following expressions 
\begin{align*}
-2\left[F^cF^c_{y}+F^sF^s_{y}\right]&= 4\sum  (\beta_{i_n}) \alpha_{i_n}\prod f_i \sin(2(\beta_{i_n})y) \\
  & \, \, +2 \sum^2 \left[ (\beta_{j_m})-(\beta_{i_n})\right]\alpha_{i_n,j_m}\prod^2 f \sin(2((\beta_{j_m})-(\beta_{i_n}))y), \\
\left[F^cF^c_{x}+F^sF^s_{x}\right] &=-2\sum  (\beta_{i_n}) \alpha_{i_n}\prod f \cos(2(\beta_{i_n})y) \\
& \, \, -2\sum^2  (\beta_{i_n}) \alpha_{i_n,j_m}\prod^2 f \cos(2((\beta_{i_n})-(\beta_{j_m}))y),
\end{align*}
\begin{align*}
\lefteqn{ -2\left[F^cF^c_{x}+F^sF^s_{x}\right]\left[F^cF^c_{y}+F^sF^s_{y}\right]= -8\sum^2  (\beta_{i_n})(\beta_{j_m}) } \\
&\qquad \quad \times \alpha_{i_n,j_m}\prod^2 f \sin(2(\beta_{i_n})y)\cos(2(\beta_{j_m})y) \\
& \qquad -4\sum^3  (\beta_{i_n}) \left[(\beta_{j_m})-(\beta_{k_p})\right] \\
& \qquad \quad  \times \alpha_{i_n,j_m,k_p}\prod^3 f\sin(2((\beta_{j_m})-(\beta_{k_p}))y)\cos(2(\beta_{i_n})y) \\ 
& \qquad -8\sum^3  (\beta_{i_n})(\beta_{j_m}) \alpha_{i_n,j_m,k_p}\prod^3 f \sin(2(\beta_{i_n})y)\cos(2((\beta_{j_m})-(\beta_{k_p}))y) \\ 
& \qquad -4\sum^4  (\beta_{i_n}) \left[(\beta_{k_p})-(\beta_{l_q})\right] \\
&\qquad \quad \times \alpha_{i_n,j_m,k_p,l_q}\prod^4 f \sin(2((\beta_{k_p})-(\beta_{l_q}))y)\cos(2((\beta_{i_n})-(\beta_{j_m}))y). 
\end{align*}
Thus, the numerator of $v$ boils down to 
\begin{align*}
& 4\sum  (\beta_{i_n})^2 \alpha_{i_n}\prod f \sin(2(\beta_{i_n})y) \\ 
&+ 2\sum^2 \left[( (\beta_{j_m})^2-(\beta_{i_n})^2)\right]\alpha_{i_n,j_m}\prod^2 f \sin(2((\beta_{j_m})-(\beta_{i_n}))y) \\
&+ 8\sum^2 (\beta_{i_n})[(\beta_{i_n})-(\beta_{j_m})] \alpha_{i_n,j_m}\prod^2 f \sin(2(\beta_{i_n})y) \cos(2(\beta_{j_m})y) \\ 
&+ 4\sum^3 [(\beta_{j_m})((\beta_{j_m})-(\beta_{i_n}))-(\beta_{k_p})((\beta_{k_p})-(\beta_{i_n}))] \\
&\quad \times \alpha_{i_n,j_m,k_p}\prod^3 f \sin(2((\beta_{j_m})-(\beta_{k_p}))y) \cos(2(\beta_{i_n})y) \\
&+ 4\sum^3 [(\beta_{i_n})((\beta_{i_n})-2(\beta_{j_m}))]\alpha_{i_n,j_m,k_p}\prod^3 f \sin(2(\beta_{i_n})y) \cos(2((\beta_{j_m})-(\beta_{k_p}))y) \\
&+ 2\sum^4 [(\beta_{j_m})((\beta_{j_m})-2(\beta_{k_p}))-(\beta_{i_n})((\beta_{i_n})-2(\beta_{k_p}))]\\
& \quad \times \alpha_{i_n,j_m,k_p,l_q}\prod^4 f \sin(2((\beta_{j_m})-(\beta_{i_n}))y) \cos(2((\beta_{k_p})-(\beta_{l_q}))y). 
\end{align*}
For sake of simplicity, assume that $0<a_n < \ldots < a_1$. From the previous expression, the behaviour of $v$ when $x\rightarrow \infty$ is governed by $f_N$,
\[
v\simeq 4  \sqrt{\dfrac{g}{h_0}}a_N f_N \sin(2\beta_N y)
\]
which is positive. 

On the other hand, the leading coefficient of $v$ when $x\rightarrow -\infty$, as the coefficient in front $(f_1\ldots f_{N})^4$ vanishes at the numerator, is given by 
\[
v\simeq -4 \sqrt{\dfrac{g}{h_0}}\alpha(1,...,N)^3\alpha(1,...,N-1)a_N (f_1\ldots f_{N-1})^4 f_N^3 \sin(2\beta_N y),
\]
which is negative.  

Let us consider the case $f_1 \simeq 1$, $f_i \ll 1, i>1$. The general case is obtained by permutations of indexes. 
The leading term of $v$ corresponds to the different powers of $f_1$,
\begin{equation}\label{s1}
v\simeq \sqrt{\dfrac{g}{h_0}}f_1(1-f_1^2)[4a_1 \sin(2\beta_1 y)]
\end{equation}
which is positive in the neighborhood of $x=s_1+c_1t$ if and only if $ x-c_1t>s_1$. The neglected terms in the previous analysis don't affect the nature of the change of sign but rather, and slightly, its position. 

In the case where $f_k \simeq 1$, $f_i \gg 1$, $i<k$ and $f_i \ll 1$, $i> k$, the leading terms of the numerator are given by
\begin{align*}
& \sqrt{\dfrac{g}{h_0}}\dfrac{4h_0^3}{3}\left( 4\alpha_{k-1}^3\alpha_{k}(f_1\ldots f_{k-1})^4f_k \sin(2\beta_ky)\left[(\beta_k)(\beta_k-(\beta_{k-1}))-(\beta_{k+1})^2\right] \right. \\
&\left. \qquad \, \, \quad +4\alpha_{k}^3\alpha_{k-1}(f_1\ldots f_{k-1})^4f_k^3 \sin(2\beta_ky)\left[-(\beta_k)^2+(\beta_{k-1})(\beta_k+(\beta_{k}))\right]\right)\\ 
=& 4a_k\sqrt{\dfrac{g}{h_0}} \Big( 1-\Big(\dfrac{\alpha(1,...,k)}{\alpha(1,...,k-1)}f_k \Big)^2 \Big) \dfrac{\alpha(1,...,k)}{\alpha(1,...,k-1)}(f_1\ldots f_{k-1})^4f_k \sin(2\beta_k y).
\end{align*}
thus, 
\[
 v\simeq 4a_k\sqrt{\dfrac{g}{h_0}}\dfrac{ \Big( 1-\Big(\frac{\alpha(1,...,k)}{\alpha(1,...,k-1)}f_k\Big)^2\Big) \frac{\alpha(1,...,k)}{\alpha(1,...,k-1)}(f_1\ldots f_{k-1})^4f_k \sin(2\beta_k y)}{\left[\left(F^c\right)^2+\left(F^s\right)^2\right]^2}. 
\]
Factoring $\alpha(1,...,k-1)^4(f_1\ldots f_{k-1})^4$ at the denominator yields,
\begin{equation}\label{estimationsigne}
 v\simeq 4a_k \sqrt{\dfrac{g}{h_0}}\left( 1-\left(\frac{\alpha(1,...,k)}{\alpha(1,...,k-1)}f_k\right)^2\right) \frac{\alpha(1,...,k)}{\alpha(1,...,k-1)}f_k\sin(2\beta_k y). 
\end{equation}
We obtain that this last expression is positive in the neighbourhood of $x=s_k+c_kt$ if and only if $ x-c_kt>s_k+\frac{1}{2\beta_k}\ln(\alpha(1,...,k)/\alpha(1,...,k-1))$. The phase shift here reproducing the phase shift in the $N$-solitons solution shown in (\ref{shift}). Here again, the neglected terms may induce a small difference between the phase shift in this analysis and the one corresponding to (\ref{shift}) but the nature of the change of sign of $v$ is preserved. 

From the previous analysis, we obtain c). Indeed, $v>0$ when $x\rightarrow \infty$ and $v$ is still positive from $x\rightarrow \infty$ to a small neighborhood of $x=s_1+c_1t$, representing the position of the crest of the first soliton. From the neighbourhood of $x=s_1+c_1t$ to the neighborhood of $x=s_2+c_2t+\frac{1}{2\beta_2}\ln(\alpha(1,2))$, only (\ref{s1}) and (\ref{estimationsigne}) play an important role in the sign of $v$. The functions $f_i$ being monotonic, only one change of sign of $v$ occurs in this region and $v$ therefore become positive for values of $x$ greater than the small neighborhood of $x=s_2+c_2t+\frac{1}{2\beta_2}\ln(\alpha(1,2))$. The analysis of (\ref{estimationsigne}) shows that there a change of sign in the neighborhood of the crest of the second soliton. The same analysis can be carried through all the crests of the $N$-solitons solution until the final crest, that is, at the left of the last soliton where $v$ is always negative. 

\endpf

\section{Numerical analysis of the particles trajectory}

\begin{figure}[!th]
\begin{center}
\begin{tabular}{c}
	\includegraphics[height=4.5cm]{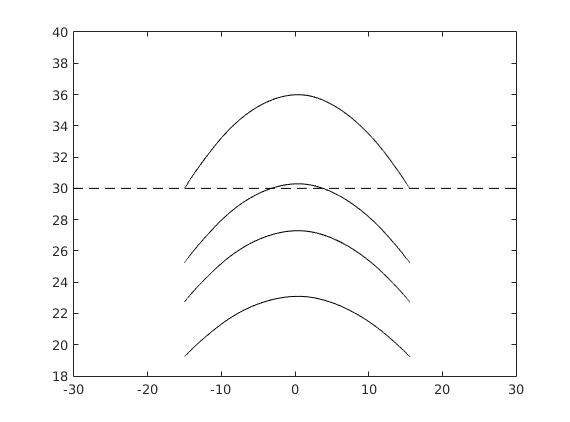}
	\includegraphics[height=4.5cm]{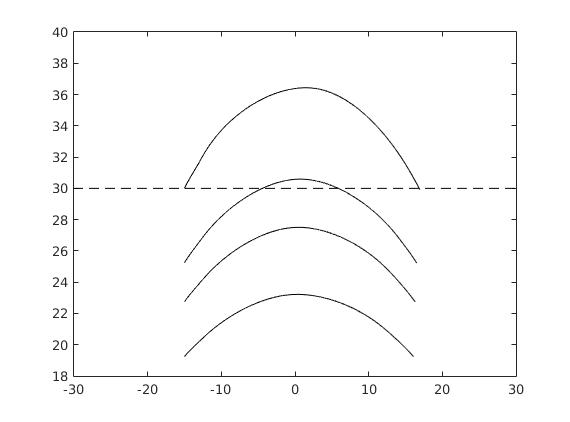} \\
	\includegraphics[height=4.5cm]{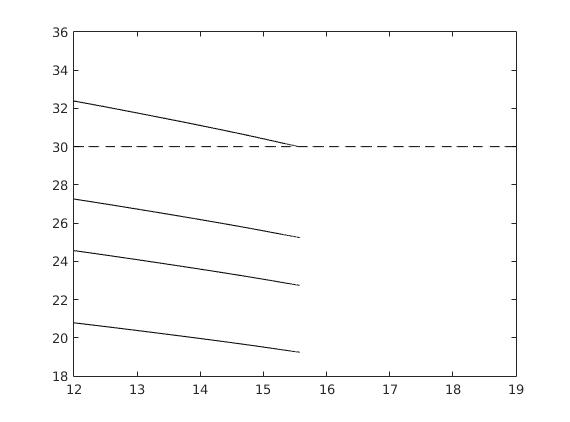}
	\includegraphics[height=4.5cm]{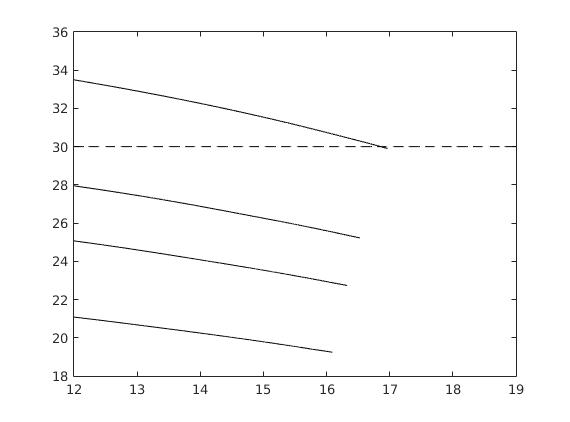}	
\end{tabular}
\end{center}
\caption{\footnotesize{Comparison of the numerical approximation of the particle trajectories for the first order velocity field (top left) and the higher order velocity field (top right). Zoom on the end of the particle trajectories for the first order velocity field (bottom left) and the higher order velocity field (bottom right). The depth of the fluid is 30 cm and the height of the solitary wave is 5.46 cm. The dashed line represents the undisturbed water surface.}}\label{c}
\end{figure}

This section is devoted to the numerical analysis of the particle trajectories given by the higher approximation of the velocity field (\ref{u})-(\ref{v}). First, we compare the results in the case of a single soliton with the numerical approximation of the particle trajectories of the first order approximation field (\ref{flow}). We use the experimental results in \cite{Chen} as the reference trajectory. The simulation were done using a fourth order Runge-Kutta scheme. Figure \ref{c} represents the numerical approximation of the particles trajectory with both velocity field with the parameters of experiment (c) of \cite{Chen}.

Table \ref{tableau} compares the total displacement in the $x$ variable and the maximal displacement in the $y$ variable of the velocity fields with the experimental results of \cite{Chen}.  
\begin{figure}
\begin{center}
    \begin{tabular}{ | l | l | l | l | l | l | l| }
   \hline
     {\footnotesize{$b$ (cm)}}& {\footnotesize{$X$ First (cm)}} &  {\footnotesize{$Y$ First (cm)}}   & {\footnotesize{$X$ Hi. (cm)}}& {\footnotesize{$Y$ Hi. (cm)}} & {\footnotesize{$X$ Exp. (cm)}} & {\footnotesize{$Y$ Exp. (cm)}} \\ \hline  
    30 & 30.57 & 6.00 & 31.97 & 6.43 & 30.63 & 5.94 \\ \hline  
    25.25 & 30.57 & 5.05 & 31.53 & 5.36 & 30.10 & 4.45 \\ \hline  
    22.75 & 30.57 &  4.55  & 31.33 & 4.77 & 30.17 & 3.93 \\ \hline  
    19.25 & 30.57 &  3.85 & 31.09 & 3.97 & 29.42 & 2.96 \\ \hline 
   \end{tabular}
   \caption{\footnotesize{Total displacement ($X$) in the x variable and maximal displacement ($Y$) in the y variable with respect to the initial vertical position above the flat bottom of the particle $b$ for the first order velocity field (First), the higher velocity field (Hi.) and the experimental results (Exp.).}}\label{tableau}
\end{center}
\end{figure} 
From the numerical results, one recovers an important feature of the particle trajectories of solitary waves : the higher the particles are initially located, the greater the horizontal displacement. However, one notices for both velocity fields that the vertical displacement of the particle located at the surface is greater that the height of the soliton. It comes from the fact that the free surface kinematic boundary conditions assuring that the particles stay underneath the water surface is not preserved in (\ref{flow}) and therefore in (\ref{u})-(\ref{v}). Indeed, this boundary condition is a higher order term in the approximation parameter as one can see from the third line of (\ref{fullvar}). Numerical simulations show that the overshoot is proportional to $\epsilon$. The numerical approximations of the particles trajectory are therefore relevant for small values of $\epsilon$ and may be interpreted qualitatively for larger values of $\epsilon$. 

\begin{figure}[!th]
\begin{center}
\begin{tabular}{c}
	\includegraphics[height=4.5cm]{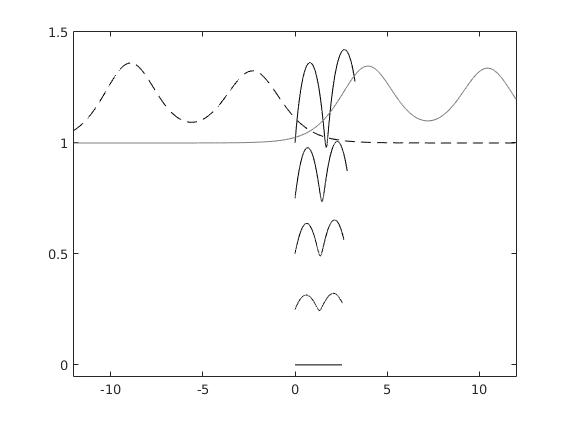}
	\includegraphics[height=4.5cm]{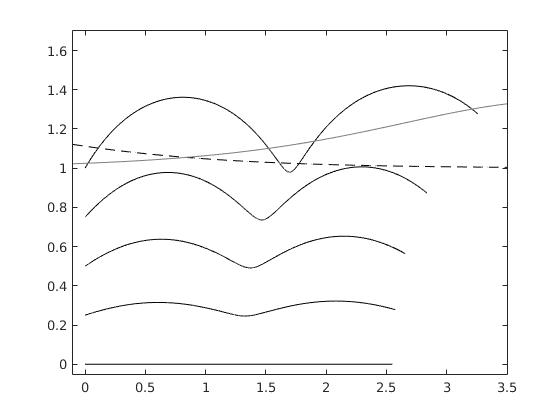}
\end{tabular}
\end{center}
\caption{\footnotesize{Numerical approximation of the particle trajectories for the 2-solitons solution. The particles trajectories are in black, the initial position of the 2-solitons is in dashed black and the final position is in gray. The height of the soliton in front is 0.4cm and the soliton behind is 0.3cm. The depth of the water is 1cm.}}\label{2sol}
\end{figure}

Figure \ref{2sol} illustrates the numerical approximation of the particles trajectory associated to the higher velocity field in the case of the 2-solitons in the case where 
\begin{equation}\label{hyp2}
 \beta_i(h_0+a)<\dfrac{\pi}{4}, \, \forall i \in \{1,\ldots,N\}
\end{equation}
holds but the sufficient condition (\ref{hyp}) of Theorem \ref{main} does not. 

Up to a renormalisation, one obtains a physically relevant particle trajectories in the case of an interaction between two solitons. Moreover, numerical investigation shows that (\ref{hyp2}) is a necessary and sufficient condition for (\ref{u})-(\ref{v}) to be physically relevant for the relevance of the particles trajectory. 

\section{Conclusion} 

It was shown theoretically that, under hypothesis (\ref{hyp}), when there is no interactions, the behaviour of the particles under the velocity field (\ref{u})-(\ref{v}) is the one induced by $N$ distinct solitons, including their phase shift. The numerical results allow us to conclude on three points which were difficult to prove theoretically for the velocity field (\ref{u})-(\ref{v}): 
\begin{enumerate}
 \item The higher the particles is initially located, the greater its horizontal displacement;
 \item Hypothesis (\ref{hyp2}), which is the theoretical bound on the amplitude parameter for solitons to exist (\cite{Whit}), is necessary and sufficient for the physical relevance of (\ref{u})-(\ref{v});
 \item The particles paths obtained from (\ref{u})-(\ref{v}) are still relevant when interactions between solitons occur.
\end{enumerate}
The numerical investigation also allowed us to highlight the possible lost of the free surface kinematic boundary condition for the first order velocity field and the higher order velocity field obtained from the first order velocity field. If one follows the derivation of the Korteweg-de Vries as in \cite{Whit}, one sees that this free surface kinematic boundary condition appears at the second order and might be recovered for second order velocity fields. However, even at this order of approximation, an overshoot is still present in some cases for second order velocity fields \cite{KalischP}. It is in contrast with the Saint-Venant equation for which the free surface kinematic boundary condition is ensured. Quantitative results are recovered for the KdV velocity fields for small value of $\epsilon$ or by renormalisation of the particles trajectories. \newline

\textbf{Acknowledgements}

The author would like to thank Jean-Michel Coron and Henrik Kalisch for the helpful discussions and the referee for the useful comments.

\bibliographystyle{plain}
\bibliography{biblio}

\end{document}